\documentclass[12pt]{article}
\usepackage{amsmath,amssymb,eucal,amsfonts,amsthm,amscd}
\newtheorem{thm}{Theorem}
\newtheorem{lem}[thm]{Lemma}
\newtheorem{prop}[thm]{Proposition}

\newtheorem{cor}[thm]{Corollary}

\newfont{\goth}{eufm10 scaled \magstep1}

\def\gg{\mbox{\goth g}}
\def\gh{{\mbox{\goth h}}}

\def\gm{\mbox{\goth m}}

\def\gt{\mbox{\goth t}}

\begin{document}

\author{Dmitri V. Alekseevsky and Liana David}

\title{A note  about invariant SKT  structures  and  generalized
K\"ahler structures on flag manifolds}

\maketitle

{\bf Abstract:} We prove that  any  invariant strong  K\"ahler
structure with torsion  (SKT  structure)  on  a flag manifold $M =
G/K$ of a  semisimple  compact Lie group $G$  is K\"ahler.
As an application we describe invariant generalized
K\"ahler structures on $M$.\\


\section{Introduction}
  A Hermitian manifold $(M,g,J)$ admits a unique  connection
  $\nabla^B$ (called the Bismut connection) which preserves  the
  metric $g$ and the complex structure $J$  and  has a skew-symmetric torsion tensor $c:= g(\cdot,
T^{B}(\cdot,\cdot))$, where
  $T^{B}$ is the torsion  of $\nabla^B$. The 3-form $c$ can be  expressed in
  terms of  the K\"ahler  form $\omega = g\circ J$  by
  $$    c = Jd\omega := d\omega(J\cdot, J \cdot,J \cdot).$$
The manifold $(M,g,J)$ is called strong K\"ahler  with
torsion (SKT) if the torsion 3-form $c$ is closed, or,
equivalently, $\partial \bar \partial \omega =0$.
SKT manifolds are  a natural generalization of K\"ahler manifolds and
many results from K\"ahler geometry can be generalized to SKT
geometry, see e.g. \cite{enrietti-fino,fernandes-fino,fino-tomassini,
fino-tomassini1}.\\

SKT  geometry is also closely related to generalized K\"ahler
geometry, which was recently introduced by N. Hitchin
\cite{hitchin} and appeared before in physics as  the geometry of
the target space of $N=(2,2)$
supersymmetric  nonlinear sigma models, see e.g. \cite{GHR,LRvUZ}.\\

A generalized K\"ahler structure on a manifold $M$
is a pair $(\mathcal{J}_1, \mathcal{J}_2)$ of commuting generalized
complex structures such that the symmetric bilinear  form
 $ -(\mathcal{J}_1 \circ \mathcal{J}_2 \cdot, \cdot)$
is positive definite, where $(\cdot , \cdot )$ is the standard
scalar product of neutral signature of the generalized tangent
bundle $\mathbb{T}M= TM\oplus T^{*}M.$ (For the definition and
basic facts about generalized complex structures, see e.g.
\cite{thesis}). It was shown by M. Gualtieri \cite{thesis,Gual}
that a generalized  K\"ahler structure on a manifold can be
described in classical terms as a bi-Hermitian  structure
$(g,J_+,J_-,b)$ in the sense of \cite{GHR}, i.e. a pair $(g,
J_+)$, $(g, J_-)$ of SKT structures with common metric $g$ and a
2-form $b$  (called in the physical literature the $b$-field) such
that
\begin{equation}\label{cond1}
db = J_+d\omega_{+}  = - J_{-}d\omega_{-},
\end{equation}
where $\omega_{\pm} = g \circ J_{\pm}$ are K\"{a}hler forms.

Let  $G$ be a semisimple compact Lie group and $M= G/K$ a
flag manifold, i.e. an adjoint orbit of $G$. In this note we describe
invariant SKT structures and invariant generalized K\"ahler structures
on $M$, as follows.

\begin{thm}\label{pmain}
 Any invariant SKT  structure $(g,J)$ on a  flag manifold $M= G/ K$ is
 K\"ahler, i.e. the K\"ahler form $\omega = g\circ J$ is closed.
\end{thm}

The description of invariant K\"ahler structures on flag
manifolds is  well known, see e.g. \cite{dmitri,dmitri'} and
Section \ref{hermitian} below.

\begin{cor}\label{mainh}  Let
$(g,J_{+}, J_{-},b)$ be an invariant  bi-Hermitian  structure
in the sense of \cite{GHR} on a flag manifold
$M = G/ K$ (which defines a generalized K\"{a}hler
structure $(\mathcal{J}_1, \mathcal{J}_2)$ via Gualtieri's
correspondence). Then  $g$ is an invariant K\"ahler metric,
$J_+,J_-$  are two parallel invariant complex structures  and
 $b$ is any closed invariant 2-form. If the group $G$ is simple,
then $J_{+}= J_{-}$ or $J_{+}= - J_{-}.$
\end{cor}

The note is organized as follows. In Section \ref{hermitian} we
fix our conventions and we recall the basic facts on the geometry
of flag manifolds and, in particular, the description of invariant
Hermitian and K\"{a}hler structures \cite{dmitri,dmitri'}. With
these preliminaries, Theorem \ref{pmain} and Corollary \ref{mainh}
will be proved in Section \ref{genkahler}.\\

{\bf Acknowledgements.} D.V.A thanks the University of Hamburg for
hospitality and financial support. L.D. acknowledges financial
support from a grant of the Romanian National Authority for
Scientific Research, CNCS-UEFISCDI, project number
PN-II-ID-PCE-2011-3-0362.

\section{Preliminary material}\label{hermitian}

{\bf Basic facts about  flag manifolds.}
 A  flag manifold of a semisimple  compact Lie  group $G$  is an
 adjoint orbit $M = \mathrm{Ad}_G (h_0) \simeq G/K$ of  an  element
 $h_0 $ of  the Lie  algebra of $G$. We denote  by
 $\mathfrak{g}$, $\mathfrak{k}$ the complex Lie algebras associated with
 the  groups $G$, $K$ respectively, and we fix
 a Cartan subalgebra $\mathfrak{h}$ of $\mathfrak{k}$. We denote by
 $R, R_0$ the  root systems  of $\mathfrak{g}$, $\mathfrak{k}$
 with respect  to $\mathfrak{h}$  and we set $R':= R \setminus R_0$.
 We write the root space decomposition of $\gg$ as
$$ \mathfrak{g} =  \mathfrak{k}+ \mathfrak{m}=( \mathfrak{h}
+ \sum_{\alpha \in R_0}\mathfrak{g}_\alpha) +
\sum_{\alpha \in R'}\mathfrak{g}_\alpha $$
 and we identify the vector space $\mathfrak{m}= \sum_{\alpha
\in R'}\mathfrak{g}_\alpha $ with the complexification of the
tangent space $T_{h_0}M $.

Let $E_{\alpha}\in \gg_{\alpha}$ be root vectors of a Weyl basis.
Thus,
$$
\langle E_{\alpha}, E_{-\alpha}\rangle =1,\quad\forall\alpha \in R
$$
(where $\langle X, Y\rangle :=\mathrm{tr}\left(\mathrm{ad}_{X}\circ
\mathrm{ad}_{Y}\right)$ denotes the Killing form of
$\gg$) and
\begin{equation}\label{weyl}
N_{-\alpha ,-\beta} = - N_{\alpha \beta},\quad\forall \alpha ,\beta
\in R
\end{equation}
where  $N_{\alpha\beta}$   are the structure  constants defined by
\begin{equation}\label{nab}
[E_{\alpha}, E_{\beta}]= N_{\alpha\beta} E_{\alpha
+\beta},\quad\forall \alpha ,\beta \in R.
\end{equation}

The  Lie algebra  of $G$  is the fixed point set
$\mathfrak{g}^\tau$
 of  the  compact anti-involution $\tau$, which
preserves the Cartan subalgebra $\gh$ and
sends $E_{\alpha}$ to  $- E_{-\alpha}$, for any $\alpha\in R.$
It is given by
 $$   \mathfrak{g}^\tau = i\gh_{\mathbb{R}} + \sum_{\alpha \in R}\mathrm{span}
 \{ E_{\alpha}-E_{-\alpha}, i(E_{\alpha}+ E_{-\alpha})\}$$
where $\gh_{\mathbb{R}} = \mathrm{span}\{ H_\alpha:=[E_{\alpha},
E_{-\alpha}],\alpha \in R\}$ is a real form of $\mathfrak{h}$.
Note that $$  \beta(H_\alpha) = \beta([E_\alpha, E_{-\alpha}]) = \langle\beta,
\alpha\rangle $$
where $\langle \cdot , \cdot \rangle$ denotes also the
scalar product on
$\mathfrak{h}^*$ induced by the Killing  form.

\subsection{Invariant complex structures on $M = G/K$}

 We  fix a system  $\Pi_0$  of simple roots of $R_0$ and we extend it
 to a system $\Pi = \Pi_0 \cup \Pi'$ of simple roots of $R$. We
 denote by $R_0^+, R^+ = R_0^+ \cup R'_+$  the corresponding
 systems of positive roots. A decomposition
 $$
\mathfrak{m} = \mathfrak{m}^+ + \mathfrak{m}^- =
 \sum_{\alpha \in R'_+}\mathfrak{g}_\alpha
+\sum_{\alpha \in R'_+}\mathfrak{g}_{-\alpha}  $$
defines an $\mathrm{Ad}_K$-invariant  complex  structure $J$ on
$T_{h_0}M = \mathfrak{m}^\tau $, such that
$$
J|_{\mathfrak{m}^{\pm}} =
\pm i \mathrm{Id}.
$$
We extend it to an invariant complex structure on $M$, also
denoted by $J$. We will refer to $R'_{+}$ and $\Pi'$ as the set of
positive roots, respectively the set of simple roots, of $J$. It
is known that any invariant complex structure on $M$ can be
obtained by  this construction \cite{dmitri,wang}.

\subsection{$T$-roots and isotropy decomposition}

Let $\mathfrak{z} = i \mathfrak{t} \subset \mathfrak{h} $ be  the
center of the  stability subalgebra $\mathfrak{k}^\tau$. The
restriction of the roots from $R' \subset \mathfrak{h}^*$ to the
subspace  $\mathfrak{t}$   are called  $T$-roots. Denote by
$$  \kappa : R' \to  R_T, \,   \alpha \mapsto \alpha|_{\mathfrak{t}}$$
the natural projection  onto the set $R_T$ of $T$-roots.  Note that
$\alpha|_{\mathfrak{t}}=0$ for any $\alpha\in R_{0}.$
Any $T$-root $\xi$  defines  an
$\mathrm{Ad}_K$-invariant subspace
$$   \mathfrak{m}_\xi := \sum_{\alpha \in R',
\kappa(\alpha) =\xi} \mathfrak{g}_\alpha $$ of the complexified
tangent space $\mathfrak{m}$ and
$$
\mathfrak{m}  = \sum_{\xi \in R_T} \mathfrak{m}_\xi
$$
is a direct sum decomposition into non-equivalent irreducible
$\mathrm{Ad}_K$-submodules.

\subsection{Invariant metrics and  Hermitian structures}

We denote by $\omega_{\alpha}\in \gg^{*}$ the $1$-forms dual to
$E_\alpha, \,\alpha \in R $, i.e.
\begin{equation}\label{label-g}
\omega_\alpha(E_\beta) = \delta_{\alpha
\beta},\quad\omega_{\alpha}\vert_{\gh}=0.
\end{equation}
Any invariant  Riemannian metric on $M$  is defined  by an
$\mathrm{Ad}_K$-invariant Euclidean metric $g$ on
$\mathfrak{m}^\tau$,  whose  complex linear extension has the form

\begin{equation}\label{adaus}
g = - \frac{1}{2}\sum_{\alpha \in R'}
g_\alpha \omega_{\alpha}\vee  \omega_{-\alpha}
\end{equation}
where $\omega_{\alpha}\vee  \omega_{-\alpha}=
\omega_{\alpha}\otimes\omega_{-\alpha}+ \omega_{-\alpha}
\otimes\omega_{\alpha}$ is the symmetric product and $g_\xi$, $\xi
\in R_T$, is a system of positive constants associated to the
$T$-roots, $g_{\xi}= g_{-\xi}$ for any $\xi \in R_{T}$ and
$g_\alpha := g_{\kappa(\alpha)}$. Note that  the restriction  of
$g$ to $\mathfrak{m}_\xi$ is proportional to the restriction of
the Killing form,  with
coefficient of proportionality $-g_\xi$.\\
Any such metric $g$ is Hermitian with respect to any invariant
complex structure $J$ and
the corresponding K\"ahler form is given by
\begin{equation}\label{omega}
\omega = - i \sum_{\alpha \in R'_+} g_\alpha
\omega_{\alpha}\wedge  \omega_{-\alpha}
\end{equation}
where $R'_{+}$ is the set of positive roots of $J$ and in our
conventions $\omega_{\alpha}\wedge\omega_{-\alpha} =
\omega_{\alpha}\otimes \omega_{-\alpha}
-\omega_{-\alpha}\otimes\omega_{\alpha}.$

\subsection{Invariant K\"ahler  structures}
Any invariant symplectic form $\omega$ on $M$
compatible  with
an invariant complex structure $J$
as above (i.e. such that
$g:= - \omega\circ J$ is
positive definite) is associated to a 1-form
 $ \sigma \in \mathfrak{t}^* $ such that  $\langle\sigma, \alpha_i
\rangle >0$ for any $\alpha_i \in \Pi'$ (the set of simple roots
of $J$). As a form on $\mathfrak{m}$, it
is given by
$$
\omega=   \omega_{\sigma}:= -i \sum_{\alpha \in R'_+} \langle
\sigma, \alpha\rangle \omega_{\alpha} \wedge \omega_{-\alpha} .$$
 The associated K\"ahler metric $g$ has
the  coefficients $g_{\alpha}=g_{\kappa(\alpha)} = \langle\sigma,
\alpha \rangle$,  which, obviously,  satisfy the  following
linearity property:
 \begin{equation} \label{linearitycondition}
   g_{\alpha + \beta} = g_{\alpha} + g_{\beta},\,\, \forall
\alpha, \beta, \alpha + \beta \in R_+'
\end{equation}
In particular, if $\Pi' = \{ \alpha_1, \cdots , \alpha_m \}$
and
$$R' \ni \alpha \equiv k_1 \alpha_1 + \cdots + k_m
\alpha_m  \, (\mathrm{mod} R_0)  $$
then
$$
g_\alpha = k_1 g_{\alpha_1} + \cdots + k_m g_{\alpha_m}.
$$

To summarize, we get:

\begin{prop}\label{linearityProp} \cite{dmitri}
An invariant Hermitian structure $(g,J)$ on $M$ is K\"ahler if and
only if the coefficients $g_{\alpha}$ associated to $g$ by
(\ref{adaus}) satisfy the linearity property: if $\alpha, \beta,
\alpha + \beta \in R'_+$, then $g_{\alpha + \beta} = g_\alpha +
g_{\beta}$. Here $R'_{+}$ is the set of positive roots of $J$.
\end{prop}

\subsection{The  formula for  the exterior derivative}

An invariant $k$-form on $M=G/K$ can  be considered as
an $\mathrm{Ad}_K$-invariant  $k$-form  $\omega$ on the
Lie algebra $\mathfrak{g}$
such that $i_{\mathfrak{k}}(\omega )=0$.
We recall the  standard  Koszul formula  for the exterior
differential $d\omega$:

\begin{equation}\label{exteriord}
( d\omega )(X_{0}, \cdots , X_{k})=\sum_{i<j}(-1)^{i+j}\omega (
[X_{i}, X_{j}],X_{1},\cdots , \widehat{X_{i}}, \cdots ,
\widehat{X_{j}}, \cdots , X_{k}),
\end{equation}
for any $X_{i}\in \gm \subset \mathfrak{g}$. In (\ref{exteriord})
the hat means that the term is omitted.

\section{Proof of our main results}\label{genkahler}

We now prove Theorem \ref{pmain} and Corollary \ref{mainh}. We preserve
the notations from the previous sections. Let $(g,J)$ be an
invariant Hermitian structure on a flag manifold $M=G/K$. Let
$g_{\alpha}= g_{k(\alpha )}$ the positive numbers associated to
$g$ and $R'_{+}$, $\Pi'$ the set of positive (respectively,
simple) roots of $J$, like before. Let  $\omega =
g \circ J$ be the K\"ahler form. To prove the theorem, we have to
check that if the form $J d\omega$ is closed, then $g_{\alpha}$
satisfy the linearity property (\ref{linearitycondition}). We
define   the sign $\epsilon_\alpha$ of a root $\alpha \in R'= R'_+
\cup(-R'_+)$ by $\epsilon_{\alpha} = \pm 1 $  if $\alpha \in \pm
R'_{+}$. Note that $\epsilon_\alpha$ depends only on
$\kappa(\alpha)$. Now  we calculate $d\omega$  and $J d\omega $ on
basic vectors, as follows:

 \begin{lem}
 \begin{enumerate}
    \item[i)]
 \begin{equation} \label{domega}
  d \omega (E_{\alpha}, E_\beta, E_\gamma) =0 \,\, {\mathrm{ if}} \,\,
   \alpha + \beta + \gamma \neq 0
 \end{equation}
and
\begin{equation} \label{domega'}  d \omega (E_{\alpha}, E_\beta, E_{-(\alpha + \beta)}) = -i
N_{\alpha \beta} (\epsilon_\alpha g_\alpha + \epsilon_{\beta}g_\beta
-\epsilon_{\alpha + \beta}g_{\alpha + \beta} ).
\end{equation}
\item[ii)]
 \begin{equation}\label{add-d}
(J d \omega )(E_{\alpha}, E_\beta, E_{-(\alpha + \beta)})=
N_{\alpha \beta} (\epsilon_\beta \epsilon_{\alpha + \beta}
g_\alpha + \epsilon_{\alpha}\epsilon_{\alpha + \beta}g_\beta -
\epsilon_{\alpha} \epsilon_{\beta}g_{\alpha + \beta} ).
\end{equation}
 \end{enumerate}
\end{lem}

\begin{proof} Relation (\ref{domega}) follows  from
(\ref{omega}) and  (\ref{exteriord}).
 Relation (\ref{domega'}) follows from
(\ref{omega}), (\ref{exteriord}) and the following property
of $N_{\alpha\beta}$ (see Chapter 5 of \cite{helgason}): \\
if
$\alpha , \beta ,\gamma\in R$ are such that $\alpha +\beta +\gamma
=0$, then
\begin{equation}\label{suma}
N_{\alpha\beta} = N_{\beta\gamma} = N_{\gamma\alpha}.
\end{equation}
Relation (\ref{add-d}) follows from (\ref{domega'}) and
$JE_{\alpha}= i\epsilon_{\alpha}E_{\alpha}$ for any $\alpha\in
R'.$

\end{proof}

\begin{lem} Suppose  that $(g, J)$ is a SKT  structure, i.e.
$d\left(  J d\omega\right) =0$. Then
\begin{equation}\label{kahler1}
N_{\alpha\beta}^{2}\left( g_{\alpha +\beta}-g_{\alpha} -
g_{\beta}\right)+ \epsilon_{\alpha -\beta}N_{\alpha,-\beta}^{2}
\left( \epsilon_{\alpha -\beta} g_{\alpha-\beta} - g_{\alpha} +
g_{\beta}\right) =0
\end{equation}
for any  $\alpha, \beta \in R_+'$, where we assume  that
$\epsilon_{\alpha -\beta} =0$ if $\alpha - \beta \notin R'$.
\end{lem}

\begin{proof} By a direct computation, we find
\begin{align*}
-\frac{1}{2}d\left( Jd\omega \right) (E_{\alpha}, E_{\beta},
E_{-\alpha}, E_{-\beta})&= N_{\alpha\beta}^{2}\left( g_{\alpha
+\beta}- g_{\alpha}
 - g_{\beta}\right)\\
&+\epsilon_{\alpha -\beta}N_{\alpha,-\beta}^{2} \left(
\epsilon_{\alpha -\beta} g_{\alpha-\beta} -g_{\alpha} +
g_{\beta}\right).
\end{align*}
This relation implies our claim.
\end{proof}

For any root
$$
R^{\prime}_{+}\ni\alpha \equiv   k_1 \alpha_1 +  \cdots +  k_{m}  \alpha_m
 \, ( \mathrm{mod}R^{+}_{0}), \,\, \alpha_i \in
\Pi',
$$
 we define the length of $\alpha$ as $\ell(\alpha) = \sum_{i=1}^{m} k_i$.
 Note that $\ell (\alpha )$ depends only on the projection
$\kappa(\alpha)$ of $\alpha$ onto $\gt^{*}.$\\

{\bf Proof of Theorem \ref{pmain}.}
By Proposition \ref{linearityProp} we have to
check  that
\begin{equation}\label{condk}
g_{\alpha +\beta} = g_{\alpha} +g_{\beta},
\end{equation}
for any $\alpha ,\beta \in R'_{+}$ such that $\alpha +\beta\in
R'_+$. We  use
induction on the length of $\gamma = \alpha + \beta \in R'_+$.
Suppose first that $\gamma = \alpha + \beta \in R'_{+}$ has  length
two.  Then  $\alpha
- \beta \notin R'$, hence $\epsilon_{\alpha -\beta}=0$.
Identity  (\ref{kahler1})
implies (\ref{condk}).\\
Suppose now that (\ref{condk}) holds for all $\gamma =
\alpha+\beta\in R'_{+}$ with $l(\gamma )\leq k$. Let $\gamma\in
R'_{+}$ with $\ell(\gamma) = k+1$  and suppose that $\gamma =
\alpha + \beta$, where $\alpha, \beta \in R'_+$. We have to show
that
\begin{equation}\label{condk0}
g_{\gamma } = g_{\alpha} +g_{\beta}.
\end{equation}
If $\alpha -\beta \notin R'$ , our previous argument shows that
(\ref{condk0}) holds. Suppose now that $\alpha -\beta\in R'$.
Without loss of generality, we may assume  that $\alpha -
\beta  \in R'_{+}.$ Then  $\alpha= (\alpha - \beta) + \beta$  is  a
decomposition of the root $\alpha $ into a sum of two roots
from $R'_{+}.$ Since $\alpha$ has length $\leq k$, our inductive
assumption implies that $g_\alpha = g_{\alpha - \beta} + g_\beta$.
Thus the second term of the identity  (\ref{kahler1}) vanishes and we
obtain (\ref{condk0}). This concludes the proof
of Theorem \ref{pmain}.\\

{\bf Proof of Corollary \ref{mainh}.} Let $(g, J_{+}, J_{-}, b)$ be
a $G$-invariant  bi-Hermitian  structure in the sense of
\cite{GHR} on a flag manifold
$M=G/K$.  Then, by Theorem \ref{pmain},  $(g,J_{\pm})$ are two
K\"ahler structures and hence the $b$-field $b$ is closed. The complex
structures $J_{\pm}$  are parallel with respect to the Levi-Civita
connection.  If the group $G$ is simple,  the K\"ahler metric $g$  is
irreducible. The endomorphism $A = J_{1}\circ J_{2}$ is symmetric with respect to $g$ and parallel. An easy argument which uses
the irreducibility of $g$  shows that
$J_{1}= J_{2}$ or $J_{1}= - J_{2}.$  This concludes the proof of Corollary  \ref{mainh}.

DMITRI V. ALEKSEEVSKY: Edinburgh University, King's Buildings,
JCMB, Mayfield Road, Edinburgh, EH9 3JZ,UK, D.Aleksee@ed.ac.uk\\

LIANA DAVID: Institute of Mathematics "Simion Stoilow" of the
Romanian Academy; Calea Grivitei nr. 21, Sector 1, Bucharest,
Romania; liana.david@imar.ro

\end{document}